\documentclass[12pt,journal,draftcls,a4paper,onecolumn]{IEEEtran}

\usepackage{amsmath,epsfig}
\usepackage{amssymb}
\usepackage{amsfonts}

\newcommand{\sg}{\boldsymbol{s}}

\newcommand{\xg}{\boldsymbol{x}}
\newcommand{\eg}{\boldsymbol{e}}
\newcommand{\yg}{\boldsymbol{y}}

\newcommand{\ggr}{\mathbf{g}}

\newcommand{\pg}{\boldsymbol{p}}

\newcommand{\alphag}{\boldsymbol{\alpha}}
\newcommand{\lambdag}{\boldsymbol{\lambda}}
\newcommand{\Mg}{\boldsymbol{M}}
\newcommand{\mg}{\boldsymbol{M}}
\newcommand{\Ng}{\boldsymbol{N}}
\newcommand{\Sigmag}{\boldsymbol{\Sigma}}
\newcommand{\Ig}{\boldsymbol{I}}

\newcommand{\ug}{\boldsymbol{u}}
\newcommand{\pgr}{\boldsymbol{p}}

\begin{document}
\onecolumn

\title{A normalized scaled gradient method to solve non-negativity and equality constrained \\linear inverse problem --  Application \\to spectral mixture analysis}

\author{\IEEEauthorblockN{C\'eline Theys \IEEEauthorrefmark{1}, Henri Lant\'eri  \IEEEauthorrefmark{1}, Nicolas Dobigeon \IEEEauthorrefmark{2},\\
    C\'edric Richard \IEEEauthorrefmark{1},  Jean-Yves Tourneret \IEEEauthorrefmark{2} and Andr\'e Ferrari \IEEEauthorrefmark{1} \vspace{1cm}}

\IEEEauthorblockA{\normalsize \IEEEauthorrefmark{1} Laboratoire
Lagrange, Universit\'e de  Nice Sophia-Antipolis, 06108 Nice
\\Email: \{celine.theys, henri.lanteri, cedric.richard,
andre.ferrari\}@unice.fr}

\IEEEauthorblockA{\normalsize \IEEEauthorrefmark{2}
IRIT/INP-ENSEEIHT, Universit\'e de Toulouse, 31071 Toulouse Cedex 7,
France \\Email: \{nicolas.dobigeon,
jean-yves.tourneret\}@enseeiht.fr}}

\maketitle

\begin{abstract}
This paper addresses the problem  of minimizing a convex cost function under non-negativity and equality constraints, with the aim of solving the linear unmixing problem encountered in hyperspectral imagery.  This problem can be formulated as a linear
regression problem whose regression coefficients (abundances)
satisfy sum-to-one and positivity cons\-traints. A normalized scaled
gradient iterative method (NSGM) is proposed for estimating the abundances
of the linear mixing model. The positivity constraint is
ensured by the  Karush Kuhn Tucker conditions whereas the
sum-to-one constraint is fulfilled by introducing normalized
variables in the algorithm. The convergence is ensured by a
one-dimensional search of the step size. Note that NSGM can be applied to any convex cost function with non negativity  and flux constraints.
In order to compare the NSGM with the well-known fully constraint least squares (FCLS) algorithm, this latter is reformulated in term of a penalized function, which reveals
its suboptimality. Si\-mu\-la\-tions on
synthetic data illustrate the performances of the proposed algorithm in comparison with other unmixing algorithms and, more particulary, demonstrate its efficiency
when compared to the popular FCLS. Finally, results on real data are given.
\end{abstract}

\section{Introduction}
Hyperspectral and multispectral imagery have received considerable
attention in the literature (see for instance \cite{Landgrebe2003,Chang2003} and references therein).
Hyperspectral and multispectral data are collected in many
spectral bands, pro\-vi\-ding accurate information regarding the
observed scene. Recent applications benefiting from
multi/hyperspectral imagery include ecosystem  monitoring
\cite{Asner2007}, crop measure \cite{Larsolle2007} and natural
disaster analysis \cite{Kokalya2007}.  Each pixel of such images is
represented by a reflectance vector, called spectrum, which contains
the measurements associated with the different spectral bands.
However, mainly due to the spatial resolution of  current
spectro-imager, the measured pixel spectrum consists of a
mixture of several spectral signatures, usually referred to as
\textit{endmembers}, that characterize the macroscopic materials
present in this pixel \cite{Keshava2002, Bioucas2012jstars}. Identifying these
endmembers and estimating their corresponding fractions, or
\emph{abundances}, in each image pixel is the core of the linear
spectral mixture analysis (LSMA). As a first approximation, it is
widely admitted that each pixel of the image can be accurately
modeled as a linear mixture of the endmembers. Following this linear
mixing model (LMM), LSMA can be addressed following two steps.
First, it is important to identify the spectral signatures
associated to the endmembers. Very popular algorithms allowing
endmember determination are N-FINDR algorithm proposed by Winter
\cite{Winter1999} and vertex component analysis (VCA) introduced in
\cite{Nascimento2005}. The second step within LSMA is the linear
unmixing of each pixel of the image. Linear unmixing consists of
estimating the abundance of each endmember contained in a
given pixel. The linear unmixing problem is challenging because the
abundances have to satisfy sum-to-one and positivity constraints.

There are mainly two kinds of approaches which can be used to
estimate  abundances that satisfy  these constraints. The first
approach is to define appropriate prior distributions for the
abundances (satisfying the sum-to-one and positivity constraints)
and estimate the unknown parameters from the resulting joint
posterior distribution following the principles of Bayesian
inference \cite{Dobigeon2008}. However, the complexity of the
parameter posterior distribution (essentially due to the
constraints inherent to abundances) requires to develop
sophisticated sampling algorithms to compute the Bayesian
estimators. These algorithms include the Gibbs sampler or the
Metropolis-within-Gibbs algorithm \cite{Robert2005}. This approach
was for instance advocated in \cite{Dobigeon2008} and provided
interesting results. The price to pay with Bayesian unmixing
algorithms is the high computational complexity resulting from the
sampling strategy.

The second approach consists of estimating the abundances by
mi\-ni\-mi\-zing an appropriate cost function such as the least squares
criterion under sum-to-one and positivity constraints. As explained
in \cite{Chang2001}, there is no analytical solution for this
optimization problem because of the linear inequalities resulting
from the positivity constraints. However, an efficient iterative
algorithm referred to as fully constrained least square (FCLS)
algorithm has been proposed in \cite{Chang2001}. The FCLS has been
applied successfully to the unmixing of hyperspectral images.
More recently, another algorithm called projected scaled gradient method (PSGM) has been proposed in
\cite{they09} where the sum-to-one constraint is ensured by a
projection at each iteration. An important advantage of the
estimators proposed in \cite{Chang2001} and \cite{they09} is their
reduced computational cost, in regard to Bayesian strategy. However, the convergence of these algorithms to the
global minimum of the cost function of interest is generally not ensured, which is their main drawback.

This paper studies a normalized split gradient method (NSGM) for
estimating the abundances involved within the LMM under positivity and sum-to-one constraints. The convergence of the NSGM is
ensured for an appropriate step size, which makes this approach very
attractive. In order to compare the proposed NSGM with the popular FCLS algorithm, we will  show that in the FCLS algorithm, the flux constraint is ensured by including a penalty term  into the quadratic data term.  The paper is organized as follows. In Section
\ref{SGM_positivity}, the normalized scaled gradient method is
developed to minimize a general criterion under positivity
constraints. The problem of taking a flux constraint (i.e., the additivity constraint)
into account is addressed in Section \ref{SGM_positivity_sum}. The resulting iterative algorithm to perform LSMA, i.e., to solve the LMM-based unmixing problem, is derived in Section \ref{sec:hyperspectral}. In Section \ref{FCLS},  explicit form  of the flux constraint is given. Simulation results conducted on synthetic and real data
are presented in Section \ref{sec:simulations}. Conclusions are
reported in Section \ref{sec:conclusion}.

\section{Scaled gradient algorithm for positivity constraints} \label{SGM_positivity}
This section studies an iterative method referred to as scaled
gradient method (SGM) to minimize any convex criterion under positivity
constraints. In other words, first the sum-to-one
constraint is not taken into account but it will be handled in the next section. Minimizing a
convex cost function $J$ under inequality constraints can be classically
achieved by introducing the Lagrange function. The Lagrange
function ${\cal L}$ associated to the linear unmixing problem with
positivity constraints can be written as
\begin{equation*}\label{55}
 {\cal L}(\alphag,\lambdag)=J(\alphag)-\lambdag^T{\ggr}(\alphag),
\end{equation*}
where $\lambdag=\left[\lambda_1 \ldots \lambda_R\right]^T$ contains
the Lagrange multipliers and $\ggr(\alphag)=\left(g(\alpha_1) \ldots g(\alpha_R)\right)^T$.  Let us note that the method proposed hereafter is valid for any
differentiable criterion $J$. Moreover if the criterion is convex and gradient Lipschitz, as in (\ref{cost}), the proposed method
converges to the global minimum of the cost function $J$. The function $g$ has to be chosen to
express the positivity constraints.
More precisely, $g$ is an increasing function that must be positive for inactive constraints
($\alpha_r>0$), and zero for active constraints $(\alpha_r=0$). The Karush Kuhn Tucker (KKT) conditions \cite{karu39,kuhn51} at the optimum
$\left(\alphag^*,\lambdag^*\right)$ express as follows
\begin{align}
  \left[ \nabla_{\alphag} {\cal L}\left(\alphag^*,\lambdag^*\right)\right]_r&=0,\qquad\forall r, \label{56}\\
g(\alpha_r^*)&\geq 0,\qquad\forall r, \label{56b} \\
\lambda_{r}^{*}&\geq 0, \qquad\forall r,  \label{56c}\\
   \lambda_{r}^{*}\;g(\alpha_r^{*})&=0, \qquad \forall r, \label{56d}
\end{align}
where $\nabla_{\alphag} {\cal L}$ is the gradient of ${\cal L}$
with respect to (w.r.t.)  $\alphag$ and the notation $[\cdot]_{r}$ is used for the $r$th
component of a vector. As a consequence, \eqref{56} leads to
\begin{equation*}\label{510}
   \lambda_{r}^{*} \frac{\partial g\left(\alpha_r^*\right)}{\partial \alpha_r}=\left[\nabla_{\alphag}
   J\left(\alphag^*\right)\right]_{r}.
\end{equation*}
Equivalently, the $r$th Lagrange multiplier can be computed as
\begin{equation*}
\lambda_{r}^{*}= \frac{\left[\nabla_{\alphag}
   J\left(\alphag^*\right)\right]_{r}}{\frac{\partial g\left(\alpha_r^*\right)}{\partial \alpha_r}},\qquad r=1,\ldots,R.
\end{equation*}
Taking into account the properties of $g$,
\eqref{56d} reduces to:
\begin{equation}
\label{511}
   \left[\nabla_{\alphag} J\left(\alphag^*\right)\right]_{r}g\left(\alpha_{r}^*\right)=0, \qquad r=1,\ldots,R.
\end{equation}
This last equation allows an iterative algorithm to be derived to
estimate $\alphag$ under positivity constraints.

Let us note that the choice for $g$ can lead to some properties on the algorithm speed. For example, taking a power function with an exponent smaller than $1$ accelerates the algorithm (see  Appendix \ref{speed}). If the descent step-size is computed to monitor the convergence of the algorithm, an obvious choice  for $g(\cdot)$ is
$g(\alpha_r)=\alpha_r$. Then, more generally, the $r$th
component of the descent direction can be
\begin{equation}
f_r(\alphag)\hat{\alpha}_{r}\left[-\nabla_{\alphag}
   J\left(\hat{\alphag}\right)\right]_{r}
\end{equation}
where $f_r(\alphag)$ is a positive function scaling the
gradient, leading to the scaled gradient method (SGM)
\begin{equation}
\alpha_{r}^{(k+1)}=\alpha_{r}^{(k)}+\gamma_r^{(k)}
f_{r}\left(\alphag^{(k)}\right)\alpha_{r}^k\left[-\nabla_{\alphag}
   J\left(\alphag^{(k)}\right)\right]_{r}, \label{SGM}
\end{equation}
where $\gamma_r^{(k)}$ is the descent step-size that must be
adjusted to ensure convergence of the algorithm.

An interesting choice for the scaling function
$f_r(\cdot)$ initially proposed in \cite{lant02} is recalled below.
The negative gradient of any convex cost function $J(\alpha)$ with a finite minimum  can always
be expressed as the difference between two positive functions:
\begin{equation} \label{split}
-\left[\nabla_{\alphag}J\left(\alphag^{(k)}\right)\right]_r=\left[U\left(\alphag^{(k)}\right)\right]_r-\left[V\left(\alphag^{(k)}\right)\right]_r.
\end{equation}
By choosing the scaling function as
\begin{equation} \label{fonctionf}
f_{r}\left(\alphag^{(k)}\right)=\frac{1}{\left[V\left(\alphag^{(k)}\right)\right]}_{r},
\end{equation}
then equation \eqref{SGM} becomes:
\begin{equation}
\alpha_{r}^{(k+1)}=\alpha_{r}^{(k)}+\gamma_r^{(k)}\alpha_{r}^k \left(\frac{\left[U\left(\alphag^{(k)}\right)\right]_r-\left[V\left(\alphag^{(k)}\right)\right]_r}{\left[V\left(\alphag^{(k)}\right)\right]_r}\right). \label{SGM2}
\end{equation}
Let us determine the maximum value for the step size in order that $\alpha_{r}^{(k+1)} \geq 0$, given $\alpha_{r}^{(k)}\geq 0$. Note that, according to \eqref{SGM2}, a restriction may only apply for the set of index $r$ such that
\begin{equation}
\label{posPQ}
    \left[U\left(\alphag^{(k)}\right)\right]_r-\left[V\left(\alphag^{(k)}\right)\right]_r<0
\end{equation}
since the other terms are positive. The maximum step size which ensures the positivity of $\alpha_{r}^{(k+1)} \geq 0$ is given by
\begin{equation}
\label{maxstep}
(\gamma_{r}^{k})_{\max}=\left[1-\frac{\left[U\left(\alphag^{(k)}\right)\right]_r}{\left[V\left(\alphag^{(k)}\right)\right]_r}\right]^{-1}
\end{equation}
which is strictly greater than $1$. Finally, the maximum step size over all the components must satisfy
\begin{equation}
\gamma^{k}_{\max}\leq \min_r \{(\gamma_{r}^{k})_{\max}\}.
\end{equation}
This choice ensures the non-negativity of the components of $\alpha_{r}^{(k)}$ from iteration to iteration. We can then write the algorithm \eqref{SGM2} with a step size $\gamma$ independent of the component
\begin{equation}
\alpha_{r}^{(k+1)}=\alpha_{r}^{(k)}+\gamma^{(k)}\alpha_{r}^k \left(\frac{\left[U\left(\alphag^{(k)}\right)\right]_r-\left[V\left(\alphag^{(k)}\right)\right]_r}{\left[V\left(\alphag^{(k)}\right)\right]_r}\right). \label{SGM3}
\end{equation}
The step size ensuring the convergence of the algorithm must be computed by an economic line search, i.e, following the Armijo rule (see appendix \ref{armi}), searched in the range $]0,\gamma^k_{\max}[$. The  step size can be chosen equal to $1$, $\gamma_r^{(k)}=1, \forall r=1,\ldots,R$, then, we obtain the classical multiplicative algorithm:
\begin{equation*}
{\alpha}_{r}^{(k+1)}=\alpha_{r}^{(k)}
\frac{\left[U\left(\alphag^{(k)}\right)\right]_r}{\left[V\left(\alphag^{(k)}\right)\right]_r}.
\label{alpha}
\end{equation*}
This multiplicative form is very attractive since the positivity of $\alpha_{r}^{(k)}$ throughout the algorithm iterations is ensured for any positive initial value $\alpha_{r}^{(0)}$ positive but the convergence is not ensured in
the general case. However, in some particular cases, for  instance
if $J\left(\hat{\alphag}\right)$ is  a quadratic cost function, we
obtain the iterative space reconstruction algorithm (ISRA)
algorithm whose convergence has been proved in \cite{daub86}.

\section{Algorithms for positivity and sum-to-one constraints} \label{SGM_positivity_sum}
\subsection{Normalized SGM}
In the case where parameters are subject to positivity and sum-to-one constraints, we propose a normalized SGM (NSGM) by
introducing the non-normalized variable $\ug=\left[u_1 \ldots
u_R\right]^T$ related to $\alphag$ by
\begin{equation}
\alphag=\frac{\ug}{\sum_j u_j}.
\end{equation}
Let us note that if $\sum_j u_j$ is a constant and if $J$ is convex
w.r.t. $\alphag$, then $J$ is also convex w.r.t.  $\ug$. This
property will be important in the following.

The gradient of $J$ w.r.t. a component $u_r$ is
\begin{equation}
\frac{\partial J}{\partial u_r}=\sum_{l=1}^R  \frac{\partial
J}{\partial \alpha_l}\frac{\partial  \alpha_l}{\partial u_r},
\end{equation}
hence
\begin{equation}
\frac{\partial J}{\partial u_r}=\frac{1}{\sum_j
u_j}\left(\frac{\partial J}{\partial \alpha_r}-\sum_{l=1}^R \alpha_l
\left( \frac{\partial J}{\partial \alpha_l}\right )\right ).
\end{equation}
Thus the negative gradient of $J$ can be decomposed as in \eqref{split} with
\begin{eqnarray}
&& \left[U\left(\ug^{(k)}\right)\right]_r =\frac{1}{\sum_j u^{(k)}_j}\left(\frac{-\partial J}{\partial \alpha_r}-\min_r\left( -\frac{\partial J}{\partial \alpha_r}\right)+\epsilon  \right ),\label{U} \\
&&\left[V\left(\ug^{(k)}\right)\right]_r =\frac{1}{\sum_j u^{(k)}_j}\left(\sum_{l=1}^R \alpha_l \left( \frac{-\partial J}{\partial \alpha_l}\right)-\min_r \left( -\frac{\partial J}{\partial \alpha_r}\right)+ \epsilon \right ).
\label{V}
\end{eqnarray}
Note that the term $\min_r\left(-\frac{\partial J}{\partial \alpha_r}\right)$
is subtracted to ensure positivity of both $U$ and $V$
whereas $\epsilon$ is a small fixed constant, that does not
modify the gradient but avoid the division by zero in
\eqref{SGM3}. Then the final algorithm based on \eqref{SGM2} with
the particular choice for $f(\cdot)$ given by
\eqref{fonctionf} and expressions \eqref{U}, \eqref{V} is:
\begin{equation}
u_{r}^{(k+1)}=u_{r}^{(k)}+\gamma_r^{(k)}u_{r}^{(k)}\left(\frac{\frac{-\partial J}{\partial \alpha_r}-\min_r\left( -\frac{\partial J}{\partial \alpha_r}\right)+\epsilon  }{\sum_{l=1}^R \alpha_l \left( \frac{-\partial J}{\partial \alpha_l}\right)-\min_r \left( -\frac{\partial J}{\partial \alpha_r}\right)+ \epsilon }-1\right).
\label{NSGMu}
\end{equation}
It can be easily found that the flux is maintained on $\ug$, i.e.,
\begin{equation}
\sum_j u_{j}^{(k+1)}=\sum_j u_{j}^{(k)}.
\end{equation}
As noticed above, the convexity of $J$ w.r.t. $\alphag$ is ensured.
This allows us to come back to the initial variables $\alphag$,
finally leading to the following updating rule of the proposed NSGM
\begin{equation}
\alpha_{r}^{(k+1)}=\alpha_{r}^{(k)}+\gamma_r^{(k)}
\alpha_{r}^{(k)} \\
\left(\frac{\frac{-\partial J}{\partial \alpha_r}-\min_r\left( -\frac{\partial J}{\partial \alpha_r}\right)+\epsilon  }{\sum_{l=1}^R \alpha_l \left( \frac{-\partial J}{\partial \alpha_l}\right)-\min_r \left( -\frac{\partial J}{\partial \alpha_r}\right)+ \epsilon }-1\right).\label{NSGM}
\end{equation}
The step sizes $\gamma_r^{(k)}$ are tuned following the
Armijo rule \cite{berts95} (see Appendix \ref{armi}), at each iteration $k$, to ensure the convergence of the algorithm.

One can be easily shown that the KKT conditions \eqref{56}, \eqref{56b}, \eqref{56c} and  \eqref{56d} are fullfilled
 at the solution:
\begin{itemize}
\item If the solution $\alpha_r^{\star}>0$, then from \eqref{SGM},  $ \left[-\nabla_{\alphag} J\left(\alphag^*\right)\right]_{r}=0$.
\item If $\alpha_r^{\star}=0$ and $ \left[\nabla_{\alphag} J\left(\alphag^*\right)\right]_{r}<0$, there is a contradiction because the term $(1+\gamma_r^{(k)} f_{r}\left(\alphag^{(k)})\right)\left[-\nabla_{\alphag} J\left(\alphag^{(k)}\right)\right]_{r}$  is greater than $1$ in the neighborhood of $\alpha_r^{\star}$ and the solution will never be reached.
\end{itemize}

\section{Application to LSMA} \label{sec:hyperspectral}

Within a widely admitted LSMA framework, the LMM assumes that a mixed pixel $\yg$ resulting from an
observation in $L$ spectral bands can be written as a linear
combination of $R$ endmember spectra $\mg_1,\ldots,\mg_R$
\begin{equation}
   \label{LMM} \yg=\Mg \alphag + \eg
\end{equation}
where $\Mg =\left(\mg_1 \ldots \mg_R\right)$ is the $L\times
R$ matrix of the endmember spectra, $\alphag=\left(\alpha_1 \ldots
\alpha_R\right)^T$ is the abundance vector to be estimated and $\eg$
is an additive noise. The linear unmixing problem considered in
this paper consists of estimating $\alphag$ under positivity and
sum-to-one constraints
\begin{equation*}
\alpha_r \geq 0, \  \forall r=1,\ldots,R \;\; \text{and} \;\;
\sum_{r=1}^R \alpha_r=1.
\end{equation*}
A standard assumption related to the LMM defined in \eqref{LMM} is
that the noise vector is distributed according to a Gaussian
distribution with zero-mean and covariance matrix $\Sigmag=\sigma^2
\Ig_L$, where $\Ig_L$ is the $L \times L$ identity matrix. 
Note that this statistical model assumes that the noise
variance is the same in all bands. This assumption has been
used extensively in the literature (see for instance
\cite{Chang1998b,Manolakis2001,Wang2006b}). Since the variance
$\sigma^2$ can be easily estimated from the observation vector, it
is assumed to be known in this paper. After removing the additive
and multiplicative constants, the resulting negative log-likelihood
function associated to the observed model reduces to the
following cost function
\begin{equation} \label{cost}
J(\alphag) =\frac{1}{2} (\yg- \Mg \alphag)^T(\yg- \Mg \alphag).
\end{equation}
The opposite of the gradient is then
\begin{equation}
\frac{\partial J}{\partial \alpha_r}=\Mg^T\yg-\Mg^T\Mg \alphag
\end{equation}
and the resulting iterative algorithm on $\alphag$ to conduct LSMA is given by the following updating rule
\begin{equation}
\alpha_{r}^{(k+1)}=\alpha_{r}^{(k)}+\gamma_r^{(k)}
\alpha_{r}^{(k)} \\
\left(\frac{\frac{-\partial J}{\partial \alpha_r}-\min_r\left( -\frac{\partial J}{\partial \alpha_r}\right)+\epsilon  }{\sum_{l=1}^R \alpha_l \left( \frac{-\partial J}{\partial \alpha_l}\right)-\min_r \left( -\frac{\partial J}{\partial \alpha_r}\right)+ \epsilon }-1\right).\label{NSGM}
\end{equation}

\section{Interpretation of the FCLS algorithm} \label{FCLS}
The FCLS algorithm \cite{Chang2001} is popular tool to solve the linear unmixing problem detailed in section \ref{sec:hyperspectral}, briefly recalled below
\begin{equation}
\min_{\alphag} (\yg- \Mg \alphag)^T(\yg- \Mg \alphag)
\end{equation}
with
\begin{equation*}
\alpha_r \geq 0, \  \forall r=1,\ldots,R \;\; \text{and} \;\;
\sum_{r=1}^R \alpha_r=1.
\end{equation*}
Within the FCLS scheme, the positivity constraint is ensured using the nonnegative least squares (NNLS) method, proposed by Lawson and Hanson \cite{laws74}. Regarding the sum-to-one constraint, Heinz and Chang introduce in \cite{Chang2001} a new signature matrix $\Ng$ and a new observation vector $\sg$ defined by
\begin{equation}
\Ng=\left[\begin{array}{c}
\delta \Mg \\
\mathbf{1}^T
\end{array}
\right] \qquad \text{and} \qquad \sg=\left[\begin{array}{c}
\delta \yg \\
1
\end{array}
\right],
\end{equation}
with $\mathbf{1}=(1 \ldots 1_R)^T$. The initial unmixing problem then becomes
\begin{equation}
\min_{\alphag} (\sg- \Ng \alphag)^T(\sg- \Ng \alphag)
\end{equation}
with
\begin{equation*}
\alpha_r \geq 0, \  \forall r=1,\ldots,R.
\end{equation*}
The negative  gradient w.r.t. $\alphag$ is then
\begin{eqnarray}\label{511}
   -\nabla_{\alphag} J\left(\alphag \right)&=&\Ng^T\sg-\Ng^T\Ng \alphag, \\
   &=&\delta^2 M^T\yg+1-\delta^2 M^T M \alphag-\sum_i \alpha_i, \\
   &=& M^T\yg-M^T M \alphag+\frac{1}{\delta^2}-\frac{1}{\delta^2}\sum_i \alpha_i.
\end{eqnarray}
Note that \eqref{511} corresponds to  the negative gradient of:
\begin{equation}
 J(\alphag)=(\yg- \Mg \alphag)^T(\yg- \Mg \alphag)+\frac{1}{2\delta^2}\left(\sum_i \alpha_i-1\right)^2.
\end{equation}
This result shows that the sum-to-one constraint is taken into account within the FCLS algorithm by adding a penalization to the data fidelity term. Consequently, the flux conservation is not ensured at each iteration and  only in an approximate way at the convergence. Moreover it depends on the regularization parameter $1/2\delta^2$ that needs to be empirically tuned. To conclude, the fundamental difference with the proposed NSGM lies in the fact that we propose an interior point method: at each iterative step, the current estimates satisfy all the constraints and we search for the best estimate among the proposed solutions.

\section{Simulation results} \label{sec:simulations}
\subsection{Synthetic data}
\begin{figure}[h!]
\begin{center}
\includegraphics[width=.5\columnwidth]{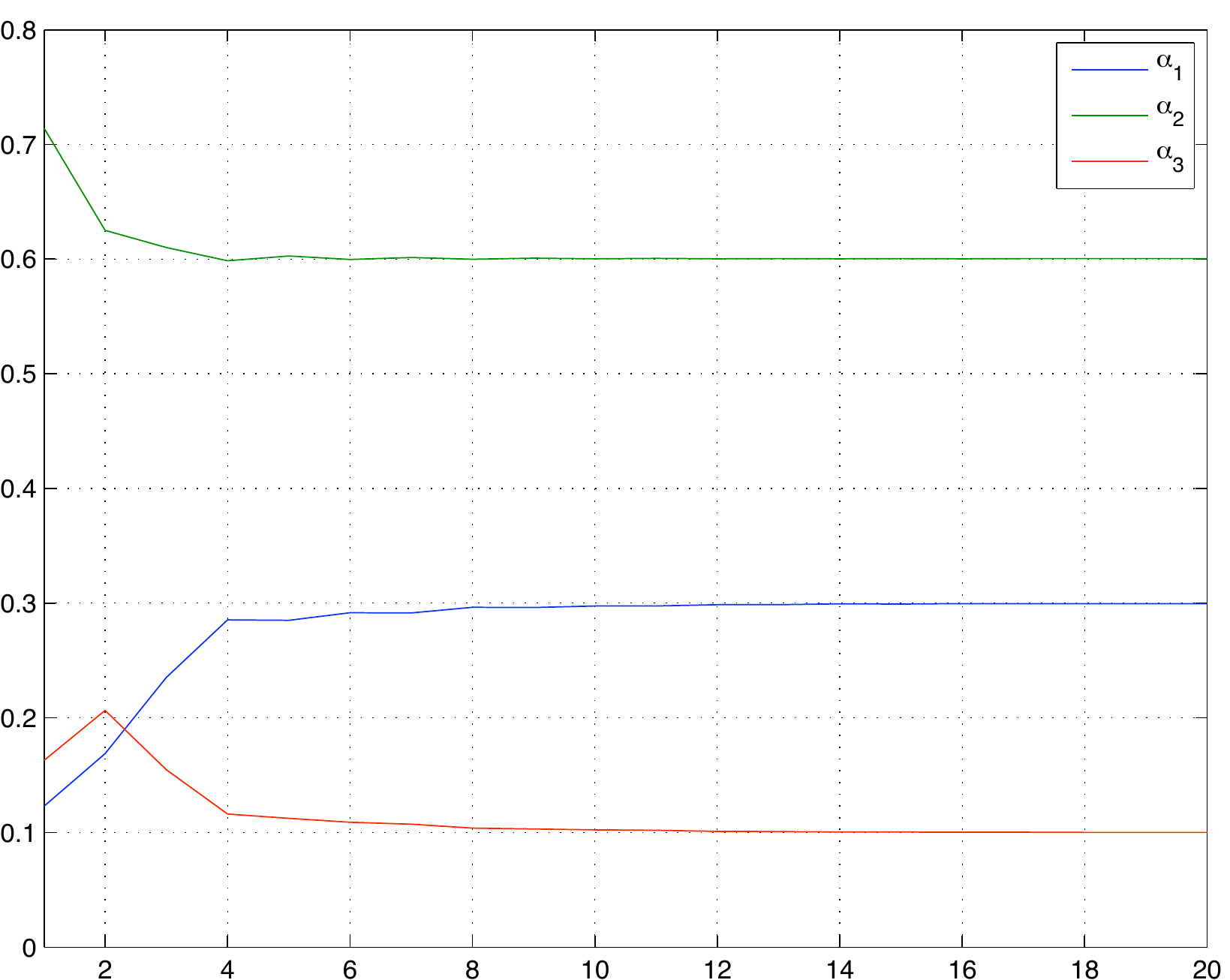}
\caption{Typical NSGM estimate $\alphag^{(k)}$ of the abundance vector
$\alphag$ versus the iteration number $k$ for $\textrm{SNR}=25$dB.
\label{NISRAk} }
\end{center}
\end{figure}

Many simulations have been conducted to validate the previous NSGM
algorithm in the LSMA context. The first experiment has been
conducted on a linear mixture of $R=3$ endmembers with
$\alphag=[0.3,0.6,0.1]^T$. The Armijo rule has been implemented with
parameters $\beta=\frac{1}{2}$ and $\sigma=\frac{1}{4}$. The three
endmembers used in this example have been extracted from the ENVI
library \cite{ENVImanual2003} and correspond to the spectra of the
construction concrete, green grass and micaceous loam. The NSGM
defined by \eqref{NSGM} has been applied on these simulated data
corrupted by an additive Gaussian noise with a signal-to-noise ratio
$\textrm{SNR}=25$dB. Figure \ref{NISRAk} shows a typical  example of
abundance estimates $\alphag^{(k)}$ as a function of the number of
iterations $k$. The algorithm clearly converges after very few
iterations.

\begin{figure}[h!]
\begin{center}
\includegraphics[width=.5\textwidth]{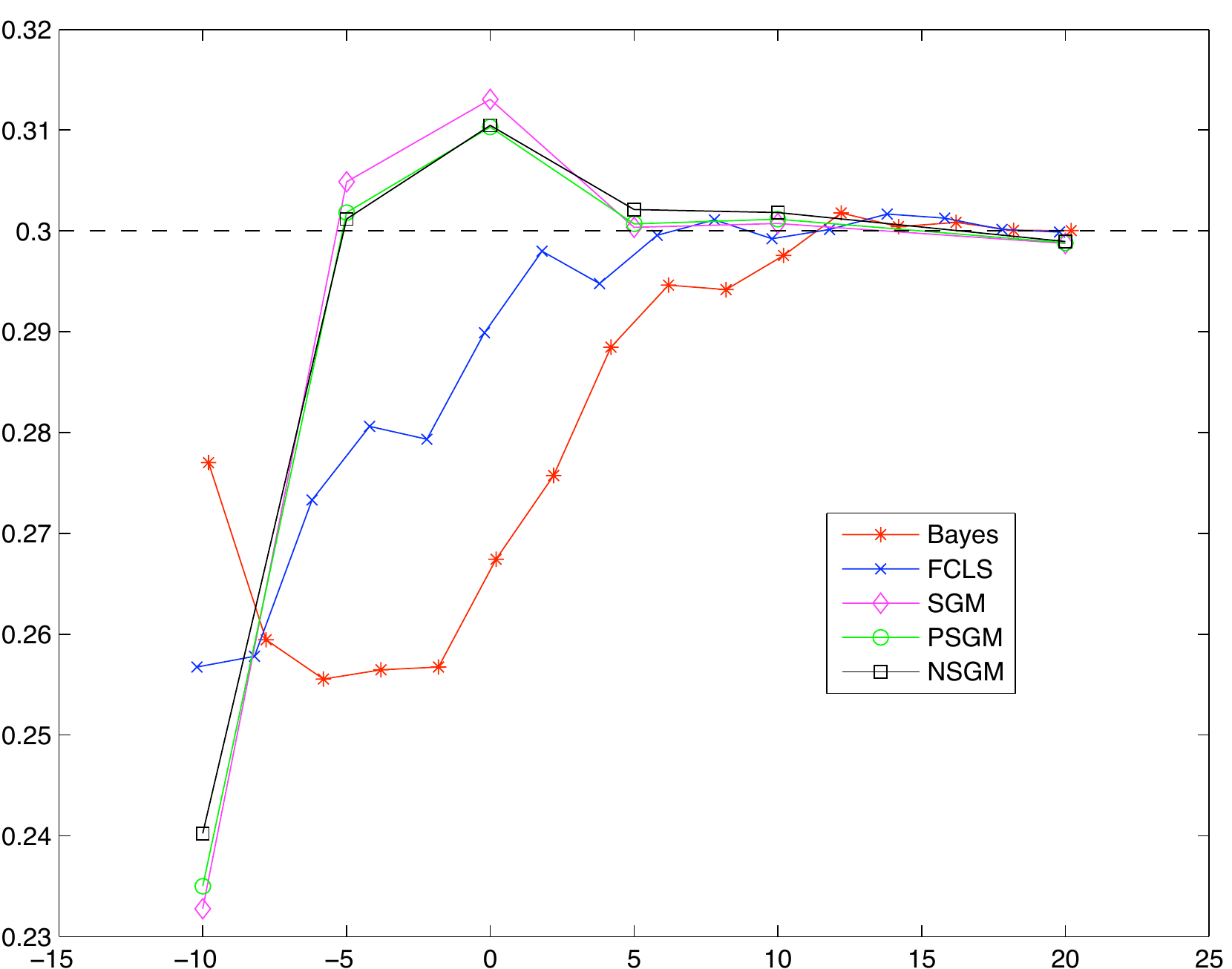}\\
\includegraphics[width=.5\textwidth]{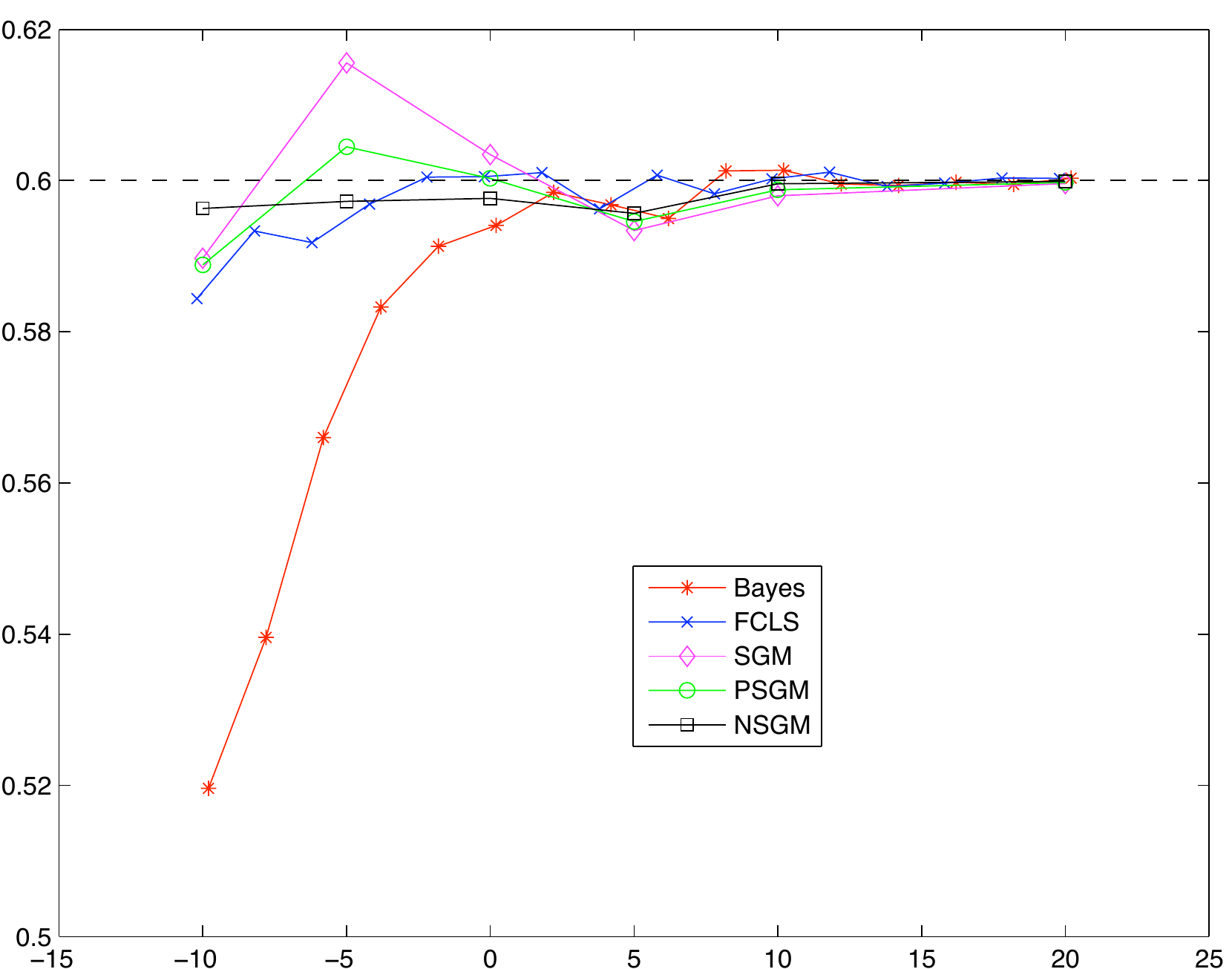}\\
\includegraphics[width=.5\textwidth]{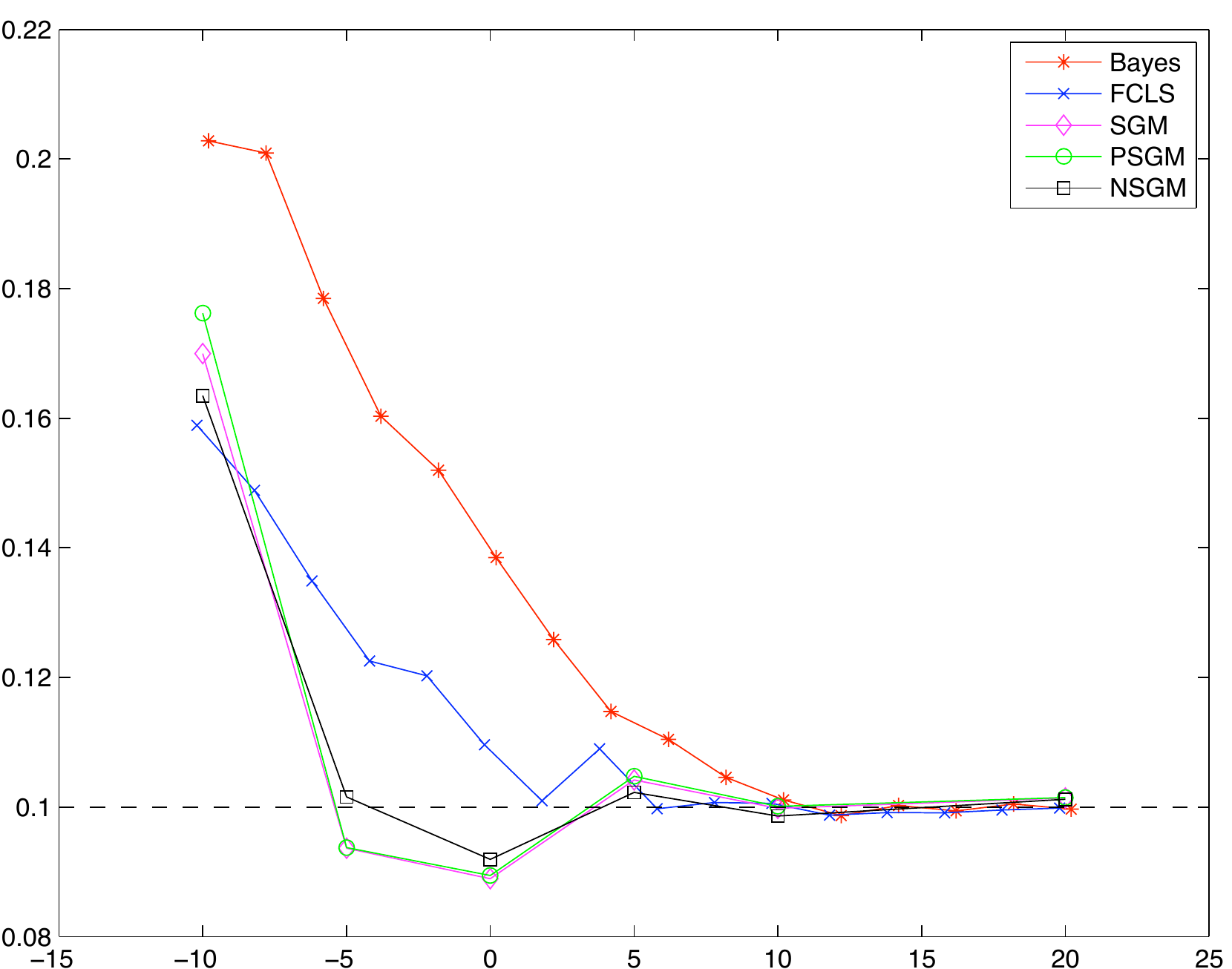}
\caption{Means  of  the abundances versus SNR for $100$ noise realizations \label{NISRA} }
\end{center}
\end{figure}

The means  of the estimated abundances with NSGM averaged over $100$ Monte Carlo runs are depicted in
Figure \ref{NISRA} as a function of the SNR and compared to those obtained with the FCLS algorithm  \cite{Chang2001}, the Bayes estimator  \cite{Dobigeon2008}, the SGM algorithm (without flux constraint) detailed in Section~\ref{SGM_positivity} and  the PSGM algorithm \cite{they09}. The corresponding estimate variances are given  in separate tables,   namely Tables \ref{tabalpha1}, \ref{tabalpha2} and \ref{tabalpha3} for parameters $\alpha_1$, $\alpha_2$ and $\alpha_3$,  respectively.

\begin{table}[h!]
\caption{Variance of $\alpha_1$ as a function of the SNR \label{tabalpha1}}
\begin{center}
{\begin{tabular}{|c|c|c|c|c|} \hline {SNR(dB)}&$-10$&$0$ &$10$&$20$
\\ \hline Bayes&$1.3e^{-2}$&$6.0e^{-3}$&$1.6e^{-3}$&$1.5e^{-4}$\\
\hline FCLS&$4.8e^{-2}$&$1.2e^{-2}$&$1.4e^{-3}$&$1.4e^{-4}$\\ \hline
SGM&$5.7e^{-2}$&$9.5e^{-3}$&$1.0e^{-3}$&$1.0e^{-4}$\\ \hline PSGM
&$5.8e^{-2}$&$8.7e^{-3}$&$1.0e^{-3}$ &$1.0e^{-4}$\\ \hline
NSGM&$5.8e^{-2}$&$8.2e^{-3}$&$1.0e^{-3}$ &$1.0e^{-4}$\\ \hline
\end{tabular}}
\end{center}
\label{default}
\end{table}

\begin{table}[h]
\caption{Variance of $\alpha_2$ as a function of the SNR \label{tabalpha2}}
\begin{center}
\begin{tabular}{|c|c|c|c|c|}
\hline {SNR(dB)}&$-10$&$0$ &$10$&$20$ \\ \hline
Bayes&$2.1e^{-2}$&$4.1e^{-3}$&$4e^{-4}$&$4e^{-5}$\\ \hline
FCLS&$3.5e^{-2}$&$4.2e^{-3}$&$4e^{-4}$&$4e^{-5}$\\ \hline
SGM&$6.6e^{-2}$&$5.1e^{-3}$&$5e^{-4}$&$5e^{-5}$\\ \hline PSGM
&$5.4e^{-2}$&$4.7e^{-3}$&$3e^{-4}$ &$4e^{-5}$\\ \hline
NSGM&$4.8e^{-2}$&$5.3e^{-3}$&$2e^{-4}$ &$3e^{-5}$\\ \hline
\end{tabular}
\end{center}
\label{default}
\end{table}%

\begin{table}[h]
\caption{Variance of $\alpha_3$ as a function of the SNR \label{tabalpha3}}
\begin{center}
\begin{tabular}{|c|c|c|c|c|}
\hline {SNR(dB)}&$-10$&$0$ &$10$&$20$ \\ \hline
Bayes&$7.7e^{-3}$&$3.3e^{-3}$&$1.1e^{-3}$&$1e^{-4}$\\ \hline
FCLS&$34.0e^{-3}$&$8.4e^{-3}$&$1.0e^{-3}$&$1e^{-4}$\\ \hline
SGM&$34.6e^{-3}$&$5.0e^{-3}$&$0.9e^{-3}$&$1e^{-4}$\\ \hline PSGM
&$37.9e^{-3}$&$5.3e^{-3}$&$0.9e^{-3}$ &$1e^{-4}$\\ \hline
NSGM&$33.0e^{-3}$&$4.5e^{-3}$&$0.8e^{-3}$ &$1e^{-4}$\\ \hline
\end{tabular}
\end{center}
\label{default}
\end{table}%

\begin{figure}[h!]
\begin{center}
\includegraphics[width=.5\textwidth]{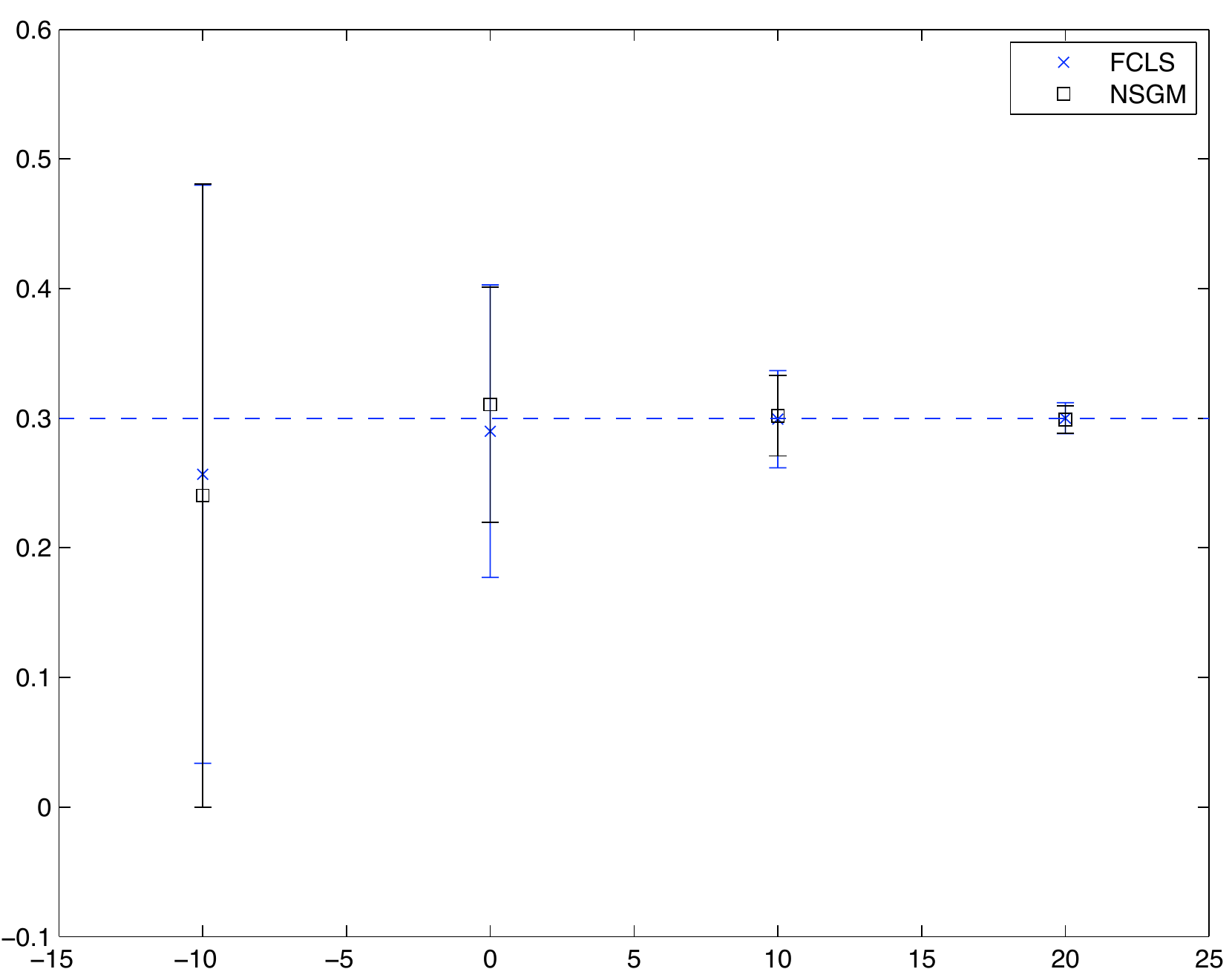}\\
\includegraphics[width=.5\textwidth]{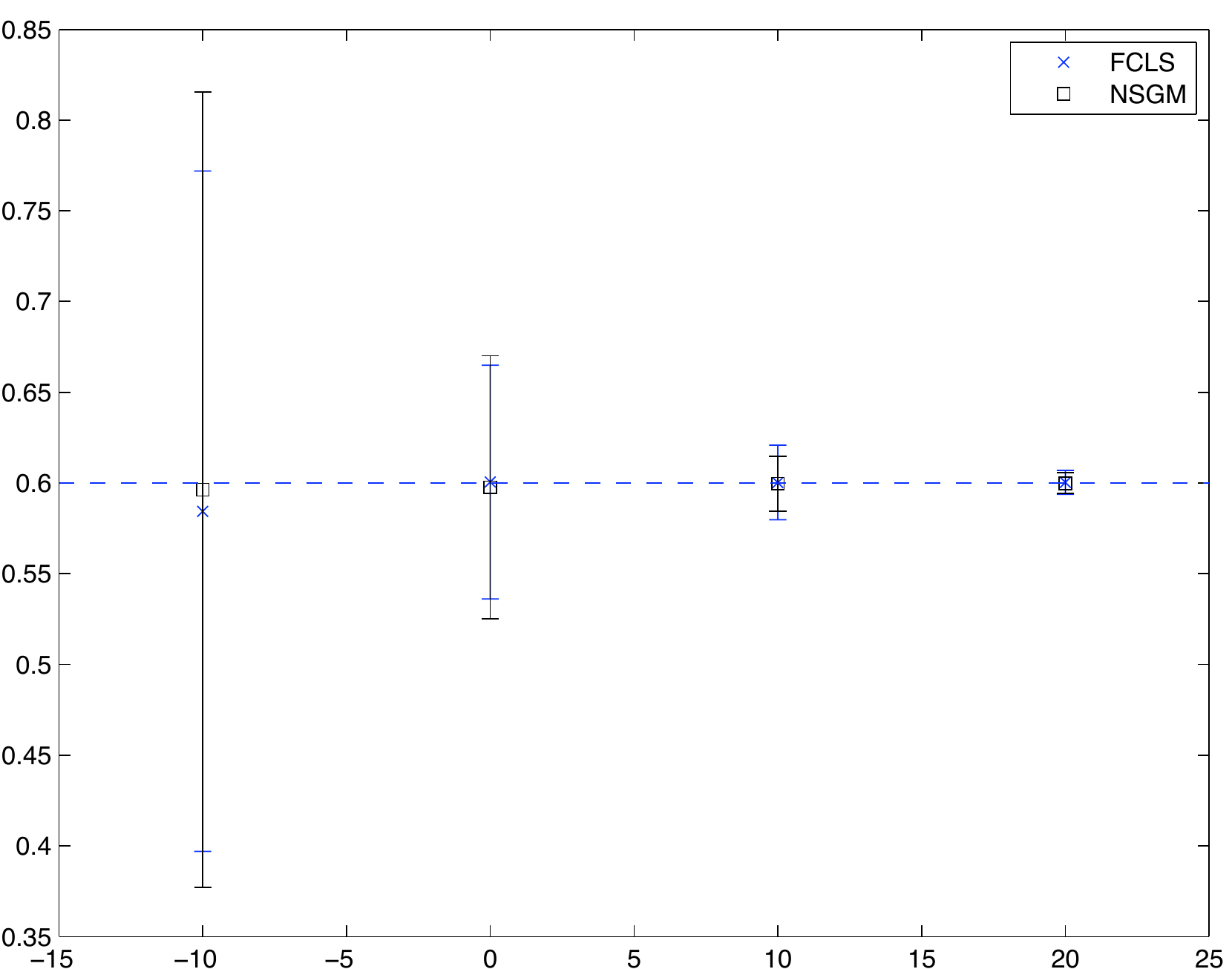}\\
\includegraphics[width=.5\textwidth]{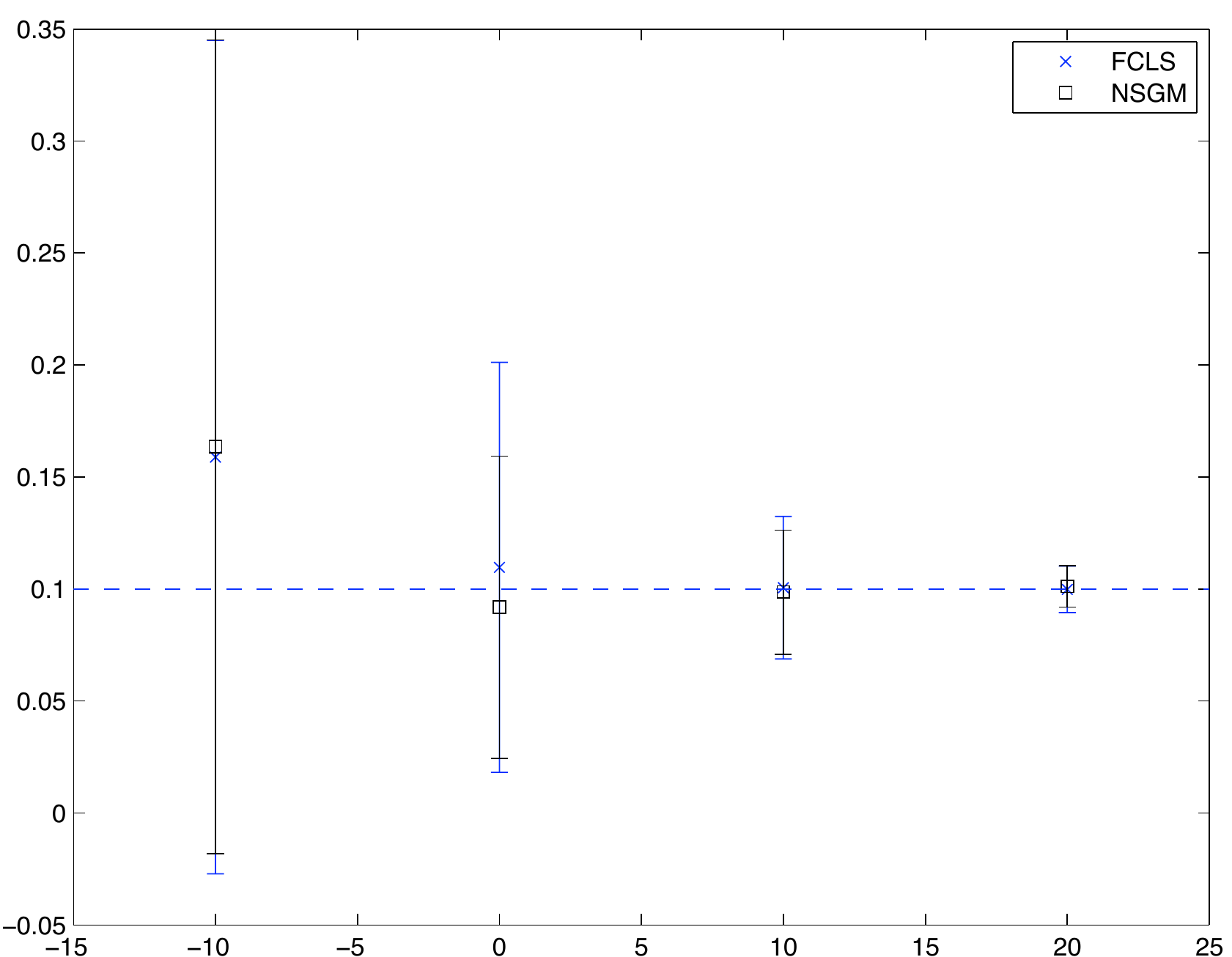}
\caption{Means  and variances of  the abundances versus SNR for FCLS and NSGM \label{FCLSNSGM} }
\end{center}
\end{figure}

Fig. \ref{FCLSNSGM} shows the means and the standard deviations of the estimated abundances for both FCLS and NSGM. Initial values of the abundances have been uniformly drawn
in the domain defined by the constraints for each Monte Carlo run. Figure  \ref{NISRA} shows that the means  obtained with the five algorithms are very similar for a SNR level over $10dB$. However, for low SNR, the iterative algorithms SGM, PSGM and NSGM perform better. Note also that their performances are very similar. However, in the case of SGM, the sum over the components of $\alphag$ fluctuates around $one$ without summing exactly to one. Moreover,  the convergence of  PSGM is not ensured contrary to the proposed NSGM.  Figure \ref{FCLSNSGM} shows that the performances of FCLS and NSGM are very similar in terms of mean and variance of the estimates but, again, the flux constraint is  strictly imposed only in the case of NSGM as it has been demonstrated in Section \ref{FCLS}.

\subsection{Real AVIRIS data}
The proposed unmixing algorithm has been also applied on two real
hyperspectral images acquired by the JPL spectroimager AVIRIS,
the Cuprite mining site (NV, USA) and the Purdue Indiana Indian test site (IN, USA).

\subsubsection{Cuprite}
\begin{figure}[h!]
\begin{center}
\includegraphics[width=1.\textwidth]{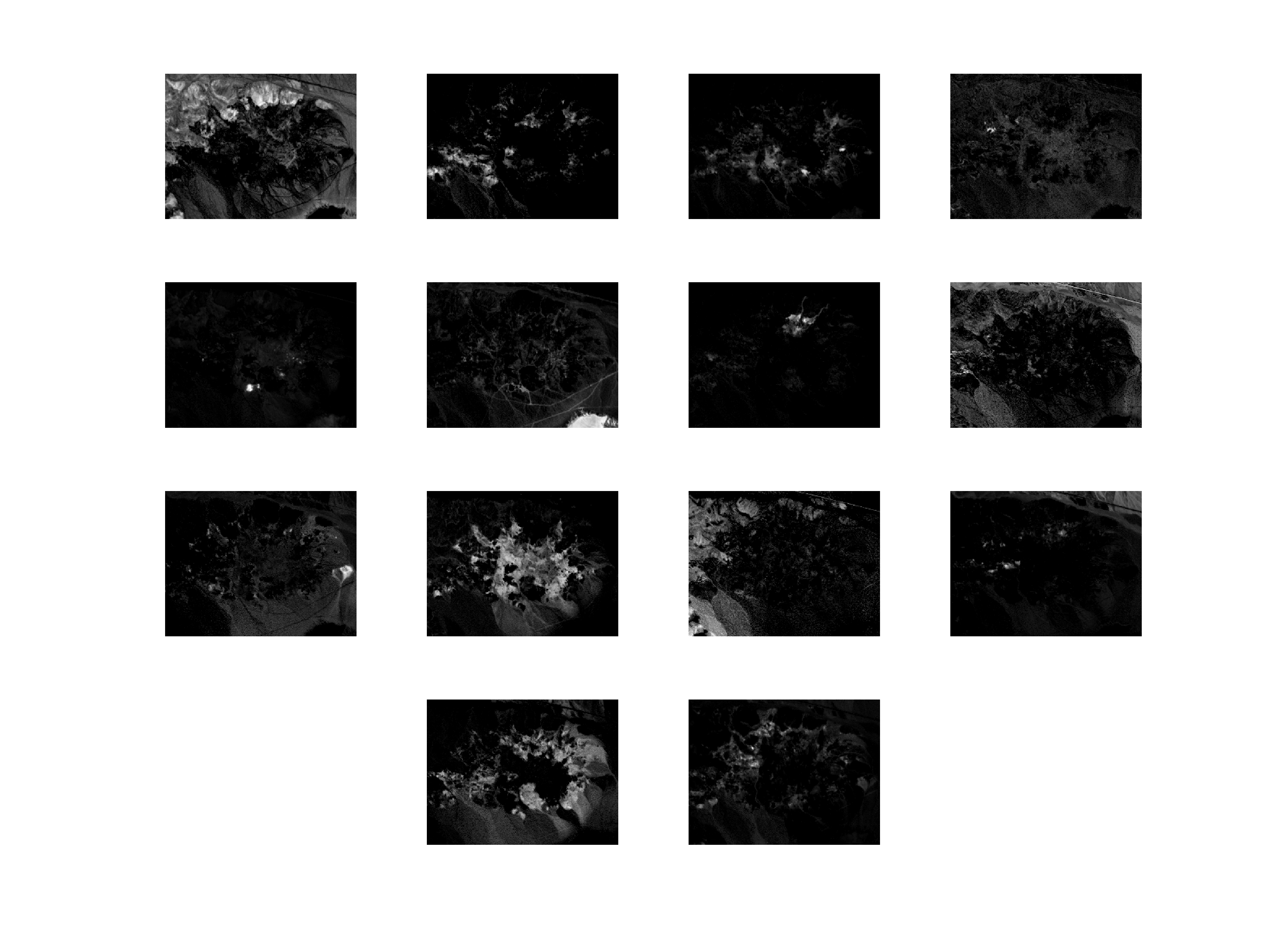}
\caption{Abundance maps estimated by the proposed normalized SGM.  \label{fig:cuprite}}
\end{center}
\end{figure}

This dataset has received
considerable attention in the literature since geologic
characteristics of the scene have been mapped in
\cite{Clark1993,Clark2003}. The sample analyzed in this experiment
consists of a sub-image of $190 \times 250$ that has been initially
studied in \cite{Nascimento2005}. Following the choice in
\cite{Nascimento2005}, $R=14$ endmember spectra have been extracted
by the vertex component analysis (VCA) \cite{Nascimento2005}.  Then, the proposed unmixing procedure has been applied
pixel-by-pixel to evaluate the relative contribution of each
endmember in each image pixel. The abundance maps are depicted in
Fig. \ref{fig:cuprite} where a black (resp. white) pixels correspond
to absence (resp. presence) of the corresponding endmembers. Note
that several areas with dominant endmembers are clearly recovered as
similar to those identified in \cite{Nascimento2005}.

\subsubsection{Indian Pines}
 The use of the abundances as features to  classify the voxels of hyperspectral image
is today a well established approach (see for example \cite{Dopido2012}). It relies
on both the dimensionality reduction and the physical meaningful of the abundance vector.

\begin{figure}[h!]
\begin{center}
\includegraphics[width=.5\textwidth]{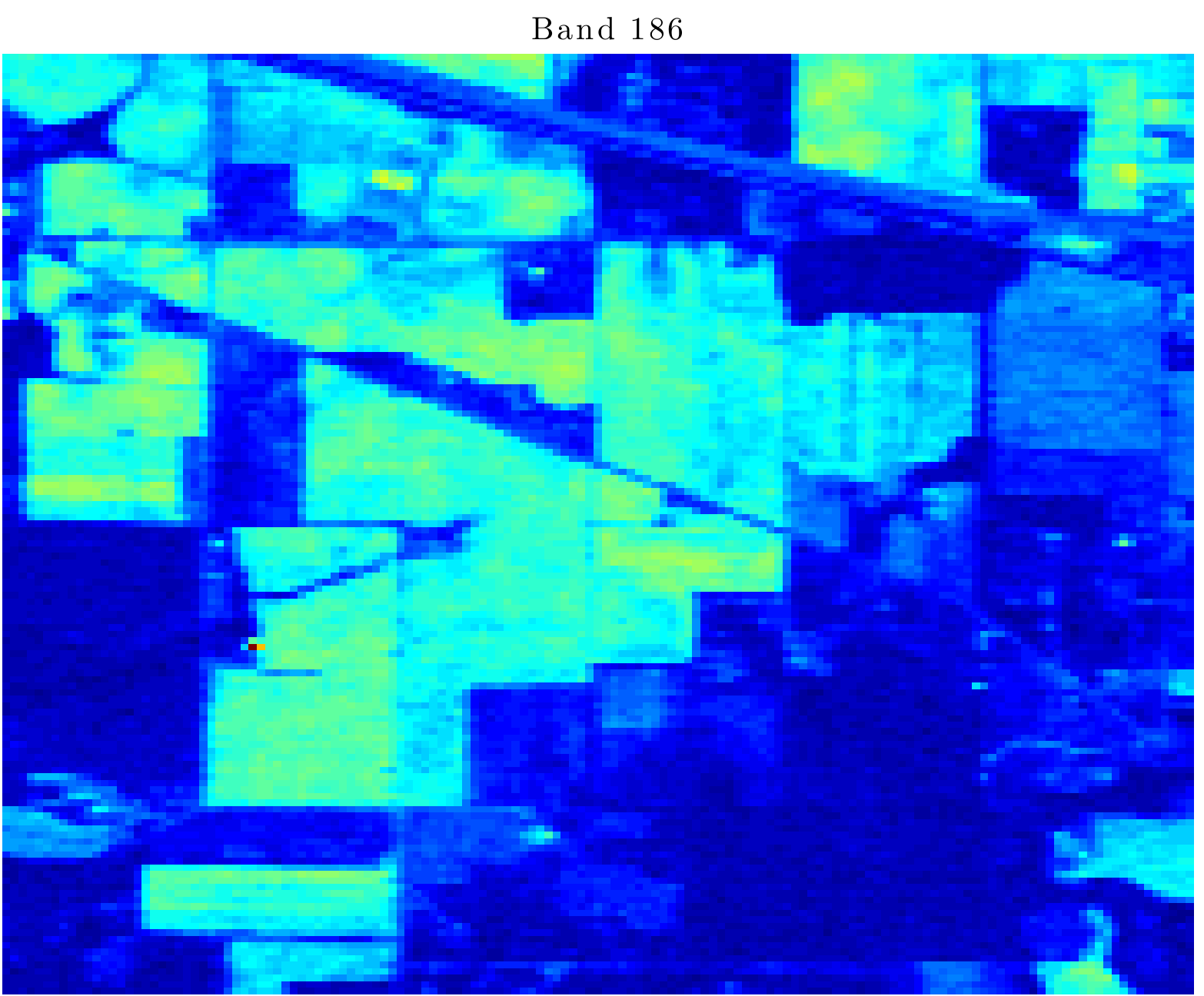}
\caption{Indian Pines data. Image at wavelength \#186.\label{band186}}
\end{center}
\end{figure}

\begin{figure}[h!]
\begin{center}
\includegraphics[width=.7\textwidth,height=8cm]{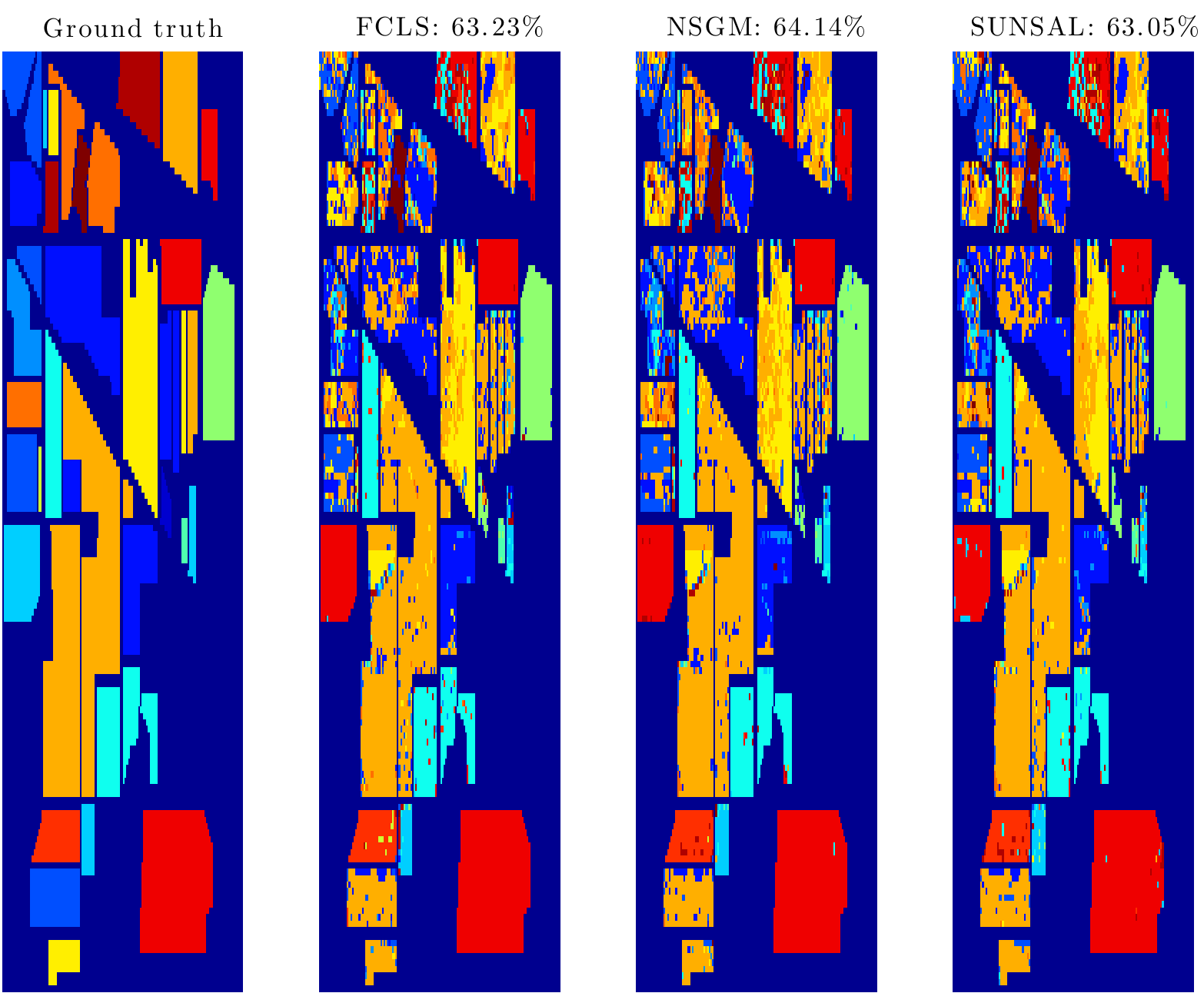}
\caption{Indian pines classification.   Fig.  (a): Ground truth. Fig. (b): Classification result using FCLS unmixing.
63.23\% of voxels are correctly classified.
Fig. (c): Classification result using NSGM unmixing. 64.14\% of voxels are correctly classified. Fig. (d):  Classification result using SUNSAL. 63.05\% of voxels are correctly classified.\label{indianpine}}
\end{center}
\end{figure}

The gain provided by the proposed unmixing algorithm in a classification framework has
been evaluated on the Indian Pines data. Figure \ref{band186} shows a typical spectral band of the $145 \times 145 \times 200$ data cube. The ground truth image shown in Figure \ref{indianpine}(a)  reveals $16$ distinct object classes and an additional background class.

The procedure followed to assess the performance of the unmixing algorithms is detailed in what follows. First, a principal component analysis has been conducted to reduce  the dimensionality of the image.  The spectra have been projected on the
subspace spanned by the $40$ principal components,
the resulting data cube containing approximately $99\%$ of the total cube power.
Then $R=18$ endmembers ($16$ corresponding to the number of classes in the ground truth, plus $2$ additional for the background class),
have been extracted using VCA. Finally, abundances have been estimated
with the FCLS algorithm \cite{Chang2001}, the SUNSAL algorithm \cite{biou2010} and the proposed fully constrained method, NSGM.

Classification of the abundances has been performed using the support-vector-machine multi-classes one-against-all
algorithm, \cite{Hsu2002}. The kernel is chosen as Gaussian with bandwidth equal to $3$ and the regularization parameter is set to $10^{-7}$.
The classifier has been trained using $10\%$ of the abundances of each class.
Figure \ref{indianpine} shows a significative result of the classification obtained by the considered methods.
Note that for this particular choice of training data $63.23\%$ of voxels are correctly classified with FCLS, $63.05\%$ with SUNSAL and $64.14\%$ with NSGM.
A Monte Carlo simulation has been performed using $100$ different training sets
randomly chosen in the image. The mean number of correctly classified voxels
is  $61.3 \%$ for FCLS, $61.7$ for SUNSAL and  $62.3\%$ for NSGM.
In this context of classification, the advantage
of the proposed method can be explained by the fact that, contrary to FCLS and SUNSAL,
NSGM estimates abundances that strictly complies to the constraints.

\section{Conclusions}\label{sec:conclusion}
Constrained scaled gradient methods were initially derived for
linear models subjected to positivity constraints.  This paper
studied  a normalized scaled gradient method (NSGM)  with  positivity and
sum-to-one constraints. NSGM can be applied to any differentiable
criterion   contrary to previous proposed algorithm, all the
constraints being fulfilled at each iteration  (characteristic to interior points methods), with an ensured convergence. The efficiency of the proposed NSGM was illustrated in a LSMA context for the
estimation of abundances. The results obtained on synthetic and real data were very promising.

\appendix

\subsection{Discussion on the relation between the  function  $g$ and the algorithm speed \label{speed}}
In Section~\ref{SGM_positivity}, the non-negativity constraint is expressed using the function $g(\alpha)=\alpha$. This appendix shows that other functions can be used, resulting in other multiplicative algorithms with higher convergence rates. Let us consider, for example, the case where the non-negativity constraint is expressed using the general function $g(\alpha)=\alpha^{1/n}$, with $n \in \cal{N^*}$. Taking the decomposition \eqref{split}:
\begin{equation}
-\left[\nabla_{\alphag}J\left(\alphag^{(k)}\right)\right]_r=\left[U\left(\alphag^{(k)}\right)\right]_r-\left[V\left(\alphag^{(k)}\right)\right]_r.
\end{equation}
The KKT condition \eqref{511} writes at the solution
\begin{equation}
\alphag_r ^{1/n} \left(\left[U\left(\alphag \right)\right]_r-\left[V\left(\alphag \right)\right]_r \right)=0
\end{equation}
that can be modified in the equivalent form
\begin{equation}
\frac{\alphag_r }{\left[V\left(\alphag \right)\right]_r^n} \left(\left[U\left(\alphag \right)\right]_r^n-\left[V\left(\alphag \right)\right]_r ^n\right).
\end{equation}
The expression
\begin{equation}
\frac{\left(\left[U\left(\alphag \right)\right]_r^n-\left[V\left(\alphag \right)\right]_r ^n\right)}{\left[V\left(\alphag \right)\right]_r^n}
\end{equation}
can be expanded in the form
\begin{equation}
\frac {\left(\left[U\left(\alphag \right)\right]_r-\left[V\left(\alphag \right)\right]_r \right)}{\left[V\left(\alphag \right)\right]_r} \left[1+\sum_{p=0}^{n-1}\left[V\left(\alphag \right)\right]_r^{p-n+1}\left[U\left(\alphag \right)\right]_r^{n-1-p}\right].
\end{equation}
Then the algorithm can be rewritten in the form
\begin{equation}
\alpha_{r}^{(k+1)}=\alpha_{r}^{(k)}+\gamma_r^{(k)}\alphag_r \frac {\left(\left[U\left(\alphag \right)\right]_r-\left[V\left(\alphag \right)\right]_r \right)}{\left[V\left(\alphag \right)\right]_r} \left[1+\sum_{p=0}^{n-1}\left[V\left(\alphag \right)\right]_r^{p-n+1}\left[U\left(\alphag \right)\right]_r^{n-1-p}\right].
\end{equation}
The function $f_r(\alphag)$ is then
\begin{equation}
f_r(\alphag)=
\frac{1}{\left[V\left(\alphag \right)\right]_r} \left[ 1+\sum_{p=0}^{n-1}\left[V\left(\alphag \right)\right]_r^{p-n+1}\left[U\left(\alphag \right)\right]_r^{n-1-p}\right].
\end{equation}

The effect of this function consists of a modification of the direction and of the modulus of the descent vector.  It is  always greater than $1/V$ and when $U$ tends to $V$, i.e, close to the convergence, it tends to $n/V$. Then taking a function $g(\alpha)=\alpha^{1/n}$ with $n>1$ has the effect to multiply the descent step-size by a factor greater than $1$ and equal to $n$ close to the convergence.  The search of the maximum step size that ensures the non-negativity of the component $\alpha$ follows the same procedure that for $n=1$ and its value is equal to
\begin{equation}
\label{maxstepn}
(\gamma_{r}^{k})_{\max}=\frac{1}{1-\frac{\left[U\left(\alphag^{(k)}\right)\right]_r^n}{\left[V\left(\alphag^{(k)}\right)\right]_r^n}}.
\end{equation}
for $r$ such that $\left[\nabla_{\alphag}J\left(\alphag \right)\right]_r >0$. It is always greater than one and in the same way,  we can obtain a multiplicative form of the algorithm  by using a constant stepsize equal to one and in this case the actualization of $\alphag$ is given by
\begin{equation*}
{\alpha}_{r}^{(k+1)}=\alpha_{r}^{(k)}
\frac{\left[U\left(\alphag^{(k)}\right)\right]_r^n}{\left[V\left(\alphag^{(k)}\right)\right]_r^n}.
\label{alpha}
\end{equation*}
Clearly, in this case, the use of the exponent $1/n$ proposed in the literature \cite{llac90,zacc96} plays the role of an accelerating term but the convergence is not ensured.
Let us note that if $0< n \leq 1$, we can easily show that
\begin{multline}
\frac{\left[U\left(\alphag^{(k)}\right)\right]^n_r-\left[V\left(\alphag^{(k)}\right) \right]^n_r}{\left[V\left(\alphag^{(k)}\right)\right]^n_r}=\frac{\left[U\left(\alphag^{(k)}\right)\right]_r-\left[V\left(\alphag^{(k)}\right) \right]_r}{\left[V\left(\alphag^{(k)}\right)\right]_r} \\
\left(\frac{1}{1+\left(\frac{\left[U\left(\alphag^{(k)}\right)\right]_r}{\left[V\left(\alphag^{(k)}\right)\right]_r}\right)^n
+\left(\frac{\left[U\left(\alphag^{(k)}\right)\right]_r}{\left[V\left(\alphag^{(k)}\right)\right]_r}\right)^{2n}+\ldots +
\left(\frac{\left[U\left(\alphag^{(k)}\right)\right]_r}{\left[V\left(\alphag^{(k)}\right)\right]_r}\right)^{(n-1)n} }\right)
\end{multline}
Then taking a function $g(\alpha)=\alpha^{1/n}$ with $0 \le n \leq 1$ has the effect to multiply the descent step size by a factor smaller than one and consequently to decrease the algorithm speed.

\subsection{Line search and Armijo rule \label{armi}}
A line search method  consists, at each iteration $k$, of choosing a descent direction $\pgr^k$ and a step length $\gamma^k$ to compute
\begin{equation}
\ug^{k+1}=\ug^k+\gamma^k \pgr^k
\end{equation}
to solve the optimization problem
\begin{equation}
\min_{\xg \in \mathbb{R}^n} f(\ug)
\end{equation}
with the following assumptions on $f(\ug)$:
\begin{itemize}
\item  $f(\ug)$ is a convex function with a finite minimum.
\item The gradient of $f(\ug)$ denoted as  $\nabla f(\ug)$, is Lipschitz continuous.
\end{itemize}

\fbox{
\begin{minipage}{\textwidth}\textbf{Armijo rule}
\begin{itemize}
\item Set scalars, $s^k,\beta,L>0,\mu$ and $\sigma$ as follows
\item \begin{equation}
s^k=\frac{-\nabla f(\ug^k)^T\pgr^k}{L||\pgr^k||^2}
\end{equation}
\item $\beta \in (0,1)$
\item $\sigma \in (0,\frac{1}{2})$
\item Then let $\gamma^k$ be the largest $\gamma$ in $\left\{s^k, \beta s^k, \beta^2 s^k, \ldots  \right\}$ such that
\begin{equation}
f(\ug^k+\gamma \pg^k)-f(\ug^k) \leq \sigma \gamma \nabla f(\ug^k)^T \pgr^k
\end{equation}
\end{itemize}
\end{minipage}
}

\newtheorem{Armijo}{Theorem}
\begin{Armijo}
Let the sequence $\{\ug\}$ be generated by $\ug^{k+1}=\ug^k+\gamma^k \pgr^k$ where $\pgr^k$ is gradient related and $\ug^k$ is chosen by Armijo rule. Then every limit point of the sequence $\left\{ \ug \right\}$ is a stationary point.
\end{Armijo}

\newpage
\bibliography{strings_all_ref,bibliosigproc}

\begin{thebibliography}{10}
\providecommand{\url}[1]{#1}
\csname url@samestyle\endcsname
\providecommand{\newblock}{\relax}
\providecommand{\bibinfo}[2]{#2}
\providecommand{\BIBentrySTDinterwordspacing}{\spaceskip=0pt\relax}
\providecommand{\BIBentryALTinterwordstretchfactor}{4}
\providecommand{\BIBentryALTinterwordspacing}{\spaceskip=\fontdimen2\font plus
\BIBentryALTinterwordstretchfactor\fontdimen3\font minus
  \fontdimen4\font\relax}
\providecommand{\BIBforeignlanguage}[2]{{%
\expandafter\ifx\csname l@#1\endcsname\relax
\typeout{** WARNING: IEEEtran.bst: No hyphenation pattern has been}%
\typeout{** loaded for the language `#1'. Using the pattern for}%
\typeout{** the default language instead.}%
\else
\language=\csname l@#1\endcsname
\fi
#2}}
\providecommand{\BIBdecl}{\relax}
\BIBdecl

\bibitem{Landgrebe2003}
D.~A. Landgrebe, \emph{Signal Theory Methods in Multispectral Remote
  Sensing}.\hskip 1em plus 0.5em minus 0.4em\relax New York: Wiley, 2003.

\bibitem{Chang2003}
C.~I. Chang, \emph{Hyperspectral Imaging: Techniques for Spectral Detection and
  Classification}.\hskip 1em plus 0.5em minus 0.4em\relax New York: Plenum
  Publishing Co., 2003.

\bibitem{Asner2007}
G.~P. Asner, D.~E. Knapp, T.~{Kennedy-Bowdoin}, M.~O. Jones, R.~E. Martin,
  J.~Boardman, and C.~B. Field, ``Carnegie airborne observatory: in-flight
  fusion of hyperspectral imaging and waveform light detection and ranging for
  three-dimensional studies of ecosystems,'' \emph{J. Appl. Remote Sensing},
  vol.~1, no.~1, p. 013536, June 2007.

\bibitem{Larsolle2007}
A.~Larsolle and H.~{Hamid Muhammed}, ``Measuring crop status using multivariate
  analysis of hyperspectral field reflectance with application to disease
  severity and plant density,'' \emph{Precision Agriculture}, vol.~8, no. 1--2,
  pp. 37--47, 2007.

\bibitem{Kokalya2007}
R.~F. Kokalya, B.~W. Rockwella, S.~L. Haireb, and T.~V.~V. Kinga,
  ``Characterization of post-fire surface cover, soils, and burn severity at
  the {C}erro {G}rande {F}ire, {N}ew {M}exico, using hyperspectral and
  multispectral remote sensing,'' \emph{Remote Sensing of Environment}, vol.
  106, no.~3, pp. 305--325, Feb. 2007.

\bibitem{Keshava2002}
N.~Keshava and J.~F. Mustard, ``Spectral unmixing,'' \emph{IEEE Signal
  Processing Magazine}, vol.~19, no.~1, pp. 44--57, Jan. 2002.

\bibitem{Bioucas2012jstars}
J.~M. Bioucas-Dias, A.~Plaza, N.~Dobigeon, M.~Parente, Q.~Du, P.~Gader, and
  J.~Chanussot, ``Hyperspectral unmixing overview: Geometrical, statistical,
  and sparse regression-based approaches,'' \emph{IEEE J. Sel. Topics Appl.
  Earth Observations and Remote Sens.}, vol.~5, no.~2, pp. 354--379, April
  2012.

\bibitem{Winter1999}
M.~E. Winter, ``Fast autonomous spectral end-member determination in
  hyperspectral data,'' in \emph{Proc. 13th Int. Conf. on Applied Geologic
  Remote Sensing}, vol.~2, Vancouver, April 1999, pp. 337--344.

\bibitem{Nascimento2005}
J.~M. Nascimento and J.~M. {Bioucas-Dias}, ``Vertex component analysis: A fast
  algorithm to unmix hyperspectral data,'' \emph{IEEE Trans. Geosci. and Remote
  Sensing}, vol.~43, no.~4, pp. 898--910, April 2005.

\bibitem{Dobigeon2008}
N.~Dobigeon, J.-Y. Tourneret, and {C.-I Chang}, ``Semi-supervised linear
  spectral unmixing using a hierarchical {B}ayesian model for hyperspectral
  imagery,'' \emph{IEEE Trans. Signal Process.}, vol.~56, no.~7, pp.
  2684--2695, July 2008.

\bibitem{Robert2005}
C.~P. Robert and G.~Casella, \emph{Monte Carlo Statistical Methods}, ser.
  Springer Texts in Statistics.\hskip 1em plus 0.5em minus 0.4em\relax New
  York: Springer-Verlag, 2005.

\bibitem{Chang2001}
D.~C. Heinz and {C.-I Chang}, ``Fully constrained least squares linear spectral
  mixture analysis method for material quantification in hyperspectral
  imagery,'' \emph{IEEE Trans. Geosci. and Remote Sensing}, vol.~39, no.~3, pp.
  529--545, March 2001.

\bibitem{they09}
C.~Theys, N.~Dobigeon, J.-Y. Tourneret, and H.~Lant{\'e}ri,
  ``\BIBforeignlanguage{anglais}{Linear unmixing of hyperspectral images using
  a scaled gradient method},'' in \emph{\BIBforeignlanguage{anglais}{Proc. IEEE
  Workshop Stat. Signal Process. (SSP)}}, Cardiff, Wales, UK, Aug. 2009, pp.
  729--732.

\bibitem{karu39}
W.~Karush, ``Minima of functions of several variables with inequalities as side
  constraints,'' Ph.D. dissertation, Univ. of Chicago, 1939.

\bibitem{kuhn51}
H.~W. Kuhn and A.~Tucker, ``Nonlinear programming,'' in \emph{Proc. 2nd
  Berkeley Symp.}, U.~of~California~Press, Ed., 1951, pp. 481--492.

\bibitem{lant02}
H.~Lant\'eri, M.~Roche, and C.~Aime, ``Penalized maximum likelihood image
  restoration with positivity constraints -- multiplicative algorithms,''
  \emph{Inverse problems}, vol.~18, no.~5, pp. 1397--1419, 2002.

\bibitem{daub86}
M.~E. Daube-Witherspoon and G.~Muehllehner, ``An iterative image space
  reconstruction algorithm suitable for volume {ECT},'' \emph{IEEE Trans.
  Medical Imaging}, vol.~5, no.~2, pp. 61--66, June 1986.

\bibitem{berts95}
D.~P. Bertsekas, \emph{Non linear programming}.\hskip 1em plus 0.5em minus
  0.4em\relax Athena Scientific, 1995.

\bibitem{Chang1998b}
{C.-I Chang}, X.-L. Zhao, M.~L.~G. Althouse, and J.~J. Pan, ``Least squares
  subspace projection approach to mixed pixel classification for hyperspectral
  images,'' \emph{IEEE Trans. Geosci. and Remote Sensing}, vol.~36, no.~3, pp.
  898--912, May 1998.

\bibitem{Manolakis2001}
D.~Manolakis, C.~Siracusa, and G.~Shaw, ``Hyperspectral subpixel target
  detection using the linear mixing model,'' \emph{IEEE Trans. Geosci. and
  Remote Sensing}, vol.~39, no.~7, pp. 1392--1409, July 2001.

\bibitem{Wang2006b}
J.~Wang and {C.-I Chang}, ``Applications of independent component analysis in
  endmember extraction and abundance quantification for hyperspectral
  imagery,'' \emph{IEEE Trans. Geosci. and Remote Sensing}, vol.~44, no.~9, pp.
  2601--2616, Sept. 2006.

\bibitem{laws74}
C.~L. Lawson and R.~J. Hanson, \emph{Solving Least Squares Problems}.\hskip 1em
  plus 0.5em minus 0.4em\relax Prentice-Hall, 1974.

\bibitem{ENVImanual2003}
{{RSI} (Research Systems Inc.)}, \emph{ENVI User's guide Version 4.0}, Boulder,
  CO 80301 USA, Sept. 2003.

\bibitem{Clark1993}
R.~N. Clark, G.~A. Swayze, and A.~Gallagher, ``Mapping minerals with imaging
  spectroscopy, {U.S. Geological Survey},'' \emph{Office of Mineral Resources
  Bulletin}, vol. 2039, pp. 141--150, 1993.

\bibitem{Clark2003}
R.~N. {{Clark \textit{et al.}}}, ``Imaging spectroscopy: {E}arth and planetary
  remote sensing with the {USGS Tetracorder} and expert systems,'' \emph{J.
  Geophys. Res.}, vol. 108, no. E12, pp. 5--1--5--44, Dec. 2003.

\bibitem{Dopido2012}
I.~D{\'o}pido, A.~Villa, A.~Plaza, and P.~Gamba, ``{A quantitative and
  comparative assessment of unmixing-based feature extraction techniques for
  hyperspectral image classification},'' \emph{IEEE J. Sel. Topics Applied
  Earth Observations and Remote Sens.}, vol.~5, no.~2, pp. 421--435, 2012.

\bibitem{biou2010}
J.~Bioucas-Dias and M.~A.~T. Figueiredo, ``Alternating direction algorithms for
  constrained sparse regression: Application to hyperspectral unmixing,'' in
  \emph{Proc. IEEE GRSS Workshop Hyperspectral Image SIgnal Process.: Evolution
  in Remote Sens. (WHISPERS)}, 2010, pp. 1--4.

\bibitem{Hsu2002}
C.-W. Hsu and C.-J. Lin, ``{A comparison of methods for multiclass support
  vector machines},'' \emph{IEEE Trans. Neur. Net.}, vol.~13, no.~2, pp.
  415--425, 2002.

\bibitem{llac90}
J.~Llacer and J.~Nu{\~n}ez, ``Iterative maximum likelihood and {B}ayesian
  algorithms for image reconstruction in astronomy,'' in \emph{The restoration
  Of Hubble Space Telescope images}, R.~L. White and R.~J. Allen, Eds.\hskip
  1em plus 0.5em minus 0.4em\relax The Space Telescope Science Institute, 1990,
  pp. 62--69.

\bibitem{zacc96}
T.~S. Zaccheo and R.~A. Gonsalves, ``Iterative maximum-likelihood estimators
  for positively constrained objects,'' \emph{J. Opt. Soc. Am. A}, vol.~13,
  no.~2, pp. 236--242, Feb 1996.

\end{thebibliography}
\bibliographystyle{ieeetran}

\end{document}